# Integrated volatility and round-off error

MATHIEU ROSENBAUM

*CREST and CMAP-École Polytechnique Paris, UMR CNRS 7641, 91128 Palaiseau Cedex, France. E-mail: mathieu.rosenbaum@polytechnique.edu*

We consider a microstructure model for a financial asset, allowing for price discreteness and for a diffusive behavior at large sampling scale. This model, introduced by Delattre and Jacod, consists in the observation at the high frequency $n$, with round-off error $\alpha_n$, of a diffusion on a finite interval. We give from this sample estimators for different forms of the integrated volatility of the asset. Our method is based on variational properties of the process associated with wavelet techniques. We prove that the accuracy of our estimation procedures is $\alpha_n \vee n^{-1/2}$. Using compensated estimators, limit theorems are obtained.

*Keywords:* diffusion models; high frequency data; integrated volatility; microstructure noise; round-off error; variation methods; wavelets

## 1. Introduction

Nowadays, a massive amount of high frequency financial data is available. This large quantity of data has paradoxically complicated some problems in statistical finance. Among them, one of the most relevant is the estimation of the integrated volatility of an asset. To fix ideas, let us consider an asset whose theoretical, efficient price $(X_t)_{t \in [0,1]}$ follows an Itô process of the form

$$dX_t = \mu_t \, dt + \sigma_t \, dW_t,$$

where $W_t$ is a Brownian motion, $\mu_t$ the drift process and $\sigma_t^2$ the instantaneous volatility. From market price observations, we wish to estimate the quantity

$$\int_0^1 g(X_t)^2 \sigma_t^2 \, dt,$$

where $g$ is a known deterministic function. The case $g(x) = 1$ corresponds to the absolute integrated volatility of the asset and the case $g(x) = 1/x$ to its relative integrated volatility, that is, the integral of the squared diffusion coefficient of the logarithmic price.[1]

---

[1]Note that the usual notion of integrated volatility refers to the relative integrated volatility.

Assume first that we observe the efficient price data with frequency $n$, that is, the sample

$$(X_{i/n}, i = 0, \ldots, n).$$

In this setting, a common convergent estimator of the integrated volatility, with rate $n^{-1/2}$ and feasible asymptotic theory, is given by the realized volatility, that is, for the absolute integrated volatility

$$\sum_{i=1}^{n}(X_{i/n} - X_{(i-1)/n})^2$$

and for the relative integrated volatility

$$\sum_{i=1}^{n}(\log(X_{i/n}) - \log(X_{(i-1)/n}))^2,$$

see Jacod and Protter [16], Barndorff-Nielsen and Shephard [7], Meddahi [22] and Gonçalves and Meddahi [12].

However, it is a well-known fact that high frequency financial data do not behave like an Itô process. In the literature, this gap is often considered to be a "contamination" of the theoretical, efficient price and is called microstructure noise. This microstructure noise increases with the sampling frequency and is due to several reasons, one of the most obvious being price discreteness.

To get rid of this noise, the first solution is to sample our data over a longer period. But, if we imagine we collect one unit of data per second and consider five minutes as the finest period we can tolerate to make the noise insignificant, we throw away a lot of data, which is hardly acceptable. Consequently, dealing with these high frequency noisy data has become a challenging issue. Many recent papers treat this problem, especially for the purpose of estimating the integrated volatility; see in particular Barndorff-Nielsen *et al.* [6]; Bandi and Russell [4]; Jacod *et al.* [15]; Zhang [30]; Zhang, Mykland and Aït-Sahalia [31]; Hansen and Lunde [13]; Aït-Sahalia, Mykland and Zhang [1]; and Gloter and Jacod [11]. For a comparison between several estimators, see Andersen, Bollerslev and Meddahi [3]; Bandi, Russel and Yang [5]; and Gatheral and Oomen [23].

In most of these works, one observes at some deterministic times $t_i^n$, $i = 0, \ldots, n$, a log-price $\tilde{Y}_{t_i^n}$ composed of a theoretical, efficient log-price $\tilde{X}_{t_i^n}$, coming from the classical continuous-time mathematical finance theory, contaminated by an additive microstructure noise $\varepsilon_{t_i^n}^n$, that is,

$$\tilde{Y}_{t_i^n} = \tilde{X}_{t_i^n} + \varepsilon_{t_i^n}^n,$$

where $\tilde{X}_t$ is, for example, an Itô process. In these additive microstructure noise models, the developed technologies often aim at reducing the impact of the noise.

Nevertheless, although price discreteness is largely accepted as one of the main reasons for microstructure noise, these models rarely allow for it; see Large [19] and Robert and Rosenbaum [25, 26] for models considering discrete prices. In this paper, we study the problem of estimating the integrated volatility of an asset when assuming that the efficient price data are observed with round-off error.

## 2. Model and results

### 2.1. Description of the model

We consider the model of a diffusion observed with round-off error. Let $\alpha_n$ be a positive decreasing sequence tending to zero as $n$ goes to infinity and $\beta_n = \alpha_n \sqrt{n}$. On a filtered probability space $(\Omega, (\mathcal{F}_t)_{t \in [0,1]}, \mathbb{P})$, we consider a one-dimensional Brownian semi-martingale $(X_t)_{t \in [0,1]}$, taking values in an open interval $(\nu, \mu)$, $-\infty \leq \nu < \mu \leq +\infty$, of the form

$$X_t = x_0 + \int_0^t \sigma(X_s)\,\mathrm{d}W_s + \int_0^t a_s\,\mathrm{d}s, \qquad (1)$$

where $(W_t)_{t \in [0,1]}$ is a $(\mathcal{F}_t)$-standard Brownian motion, $(a_t)_{t \in [0,1]}$ is a progressively measurable process with respect to $(\mathcal{F}_t)_{t \in [0,1]}$, $x \to \sigma(x)$ is a real deterministic function that is not known and $x_0$ is a real constant. We observe the sample

$$(X_{i/n}^{(\alpha_n)}, i = 0, \ldots, n), \qquad (2)$$

where

$$X_{i/n}^{(\alpha_n)} = \alpha_n \lfloor X_{i/n}/\alpha_n \rfloor.$$

Thus, $X_{i/n}^{(\alpha_n)}$ is the observation of $X_{i/n}$ with round-off error $\alpha_n$. This model has already been studied by Delattre and Jacod [9] when $\beta_n$ tends to a constant finite value and by Delattre [8] in the other cases. Based on the sample (2), our goal is to estimate the random parameter

$$\theta = \int_0^1 g(X_s)^2 \sigma(X_s)^2\,\mathrm{d}s,$$

where $g$ is a known deterministic function on $(\nu, \mu)$.

Note that for the Black–Scholes specification of the model

$$\sigma(x) = \sigma x,$$

the problem of the estimation of the constant parameter $\sigma$ has been partially treated by Li and Mykland [20] in the case where $\beta_n$ tends to zero.

We denote by $\mathcal{C}^k(I)$ the set of $k$ times continuously differentiable functions on $I \subseteq \mathbb{R}$. We write $\mathcal{C}_b^k(I)$ if all the derivatives are bounded. We will consider the following assumptions:

**Assumption A.**
$$\sup_{n \geq 0} \alpha_n (\log n)^2 < \infty.$$

**Assumption A1.** There exists $\rho > 0$ such that $\sup_{n \geq 0} \alpha_n^{1-\rho} (\log n)^2 < \infty$.

**Assumption B.**

(i) For all $x \in (\nu, \mu)$, $\sigma(x) > 0$,
(ii) $x \to \sigma(x) \in \mathcal{C}^2((\nu, \mu))$,
(iii) $\int_0^1 a_s^2 \, ds < +\infty$, almost surely.

**Assumption C.**
$$x \to g(x) \in \mathcal{C}^2((\nu, \mu)) \text{ and for all } x \in (\nu, \mu), \, g(x) > 0.$$

**Assumption C1.**

(i) $x \to g(x) \in \mathcal{C}^2((\nu, \mu))$ and for all $x \in (\nu, \mu)$, $g(x) > 0$,
(ii) $x \to g'(x)$ is of constant sign on $(\nu, \mu)$ and $x \to |g'(x)|^{1/2} \in \mathcal{C}^2((\nu, \mu))$.

By convention, if $x \notin (\nu, \mu)$, we set $g(x) = g'(x) = 0$.

Note that on $(0, +\infty)$, the functions defined by $x \to 1$ and $x \to 1/x$ satisfy Assumption C1. These functions are those respectively associated to the absolute integrated volatility and to the relative integrated volatility.

## 2.2. First estimator

Our estimation method is based on the theory of wavelet methods for quadratic functionals estimation; see, for example, Gayraud and Tribouley [10]. Throughout the paper, for $k \in \mathbb{N}$ and $j \in \mathbb{N}$, we set

$$\mathbb{1}_{jk}(s) = \mathbb{1}_{(k/2^j, k+1/(2^j)]}(s), \qquad \psi(s) = -\mathbb{1}_{[0,1/2]}(s) + \mathbb{1}_{(1/2,1]}(s),$$
$$\psi_{jk}(s) = 2^{j/2} \psi(2^j s - k).$$

We define the coefficients $c_{j_0 k}$, $j_0 \in \mathbb{N}$, $k \in [0, 2^{j_0} - 1]$ and $d_{jk}$, $j \in \mathbb{N}$, $k \in [0, 2^j - 1]$ by

$$c_{j_0 k} = 2^{j_0/2} \int \mathbb{1}_{j_0 k}(s) g(X_s) \sigma(X_s) \, ds, \qquad d_{jk} = \int \psi_{jk}(s) g(X_s) \sigma(X_s) \, ds.$$

Hence, the $c_{j_0 k}$ and $d_{jk}$ are the coefficients of $s \to g(X_s)\sigma(X_s)$ in the Haar basis. Consequently, we have

$$\theta = \sum_{k=0}^{2^{j_0}-1} c_{j_0 k}^2 + \sum_{j=j_0}^{+\infty} \sum_{k=0}^{2^j-1} d_{jk}^2.$$

We set

$$\hat{c}_{j_0 k} = \sqrt{\frac{\pi}{2}} \frac{2^{j_0/2}}{\sqrt{n}} \sum_{i=1}^{n} \mathbb{1}_{j_0 k}(i/n) g(X_{(i-1)/n}^{(\alpha_n)}) |X_{i/n}^{(\alpha_n)} - X_{(i-1)/n}^{(\alpha_n)}|.$$

Thus, in the case $g(x) = 1$, $\hat{c}_{j_0 k}$ can be seen as a rescaled local average of the increments of the rounded diffusion over a window of size $2^{-j_0}$. We define our first estimator $\widetilde{\theta}_n$ of $\theta$ by

$$\widetilde{\theta}_n = \sum_{k=0}^{2^{j_{0,n}}-1} \hat{c}_{j_{0,n} k}^2$$

with $j_{0,n} = \lfloor \log_2(\alpha_n^{-1} \wedge \sqrt{n}) \rfloor$.

### 2.3. Convergence in probability

We set $r_n = \alpha_n \vee n^{-1/2}$. We have the following theorem:

**Theorem 1 (Convergence in probability).** *In model (1)–(2), under Assumptions A, B and C, the sequence*

$$r_n^{-1}(\widetilde{\theta}_n - \theta)$$

*is tight.*

### 2.4. Compensated estimators

It seems difficult to obtain a central limit theorem for the previous estimator (see the proofs for details). Consequently, we introduce compensated estimators. We set

$$Q_j = \sum_{k=0}^{2^j - 1} d_{jk}^2$$

and define

$$\hat{d}_{jk} = \sqrt{\frac{\pi}{2}} \frac{1}{\sqrt{n}} \sum_{i=1}^{n} \psi_{jk}(i/n) g(X_{(i-1)/n}^{(\alpha_n)}) |X_{i/n}^{(\alpha_n)} - X_{(i-1)/n}^{(\alpha_n)}|.$$

We denote by $\mathcal{S}$ the set of all triples $(a, (j_{1,n}), (j_{2,n}))$, where $a$ is a real number with $0 < a < 1$ and $(j_{1,n}), (j_{2,n})$ are two sequences of integers such that

$$\sup_n \alpha_n^{1-a}(\log n)^2 < \infty, \qquad r_n 2^{2j_{2,n} - j_{1,n}} \to 0, \qquad r_n^{-1} 2^{j_{1,n}/2}(\alpha_n^2 \log n + 1/n) \to 0,$$

$$r_n 2^{j_{1,n}} \to 0, \qquad r_n^{-1} 2^{-3j_{1,n}/2} \to 0, \qquad 2^{j_{2,n} - j_{1,n}} \to 0, \qquad r_n^{-1} 2^{-(j_{1,n} + j_{2,n}/2)} \to 0.$$

Under Assumption A1, the set $\mathcal{S}$ is not empty. For example, if one takes $j_{1,n} = \lfloor \log_2(r_n^{-3/4}) \rfloor$ and $j_{2,n} = \lfloor \log_2(r_n^{-2/3}) \rfloor$, then $(\rho, (j_{1,n}), (j_{2,n})) \in \mathcal{S}$. For $S = (a, (j_{1,n}), (j_{2,n})) \in \mathcal{S}$, we set

$$\hat{Q}_{j_{2,n}} = \sum_k \hat{d}_{j_{2,n}k}^2$$

and consider

$$R_n(S) = \sum_{j=j_{1,n}}^{\lfloor (1+a)\log_2 r_n^{-1} \rfloor} 2^{j_{2,n}-j} \hat{Q}_{j_{2,n}}.$$

For $S = (a, (j_{1,n}), (j_{2,n})) \in \mathcal{S}$, our final estimator of $\theta$ is

$$\hat{\theta}_n(S) = \sum_{k=0}^{2^{j_{1,n}}-1} \hat{c}_{j_{1,n}k}^2 + R_n(S) + \alpha_n(\mathbb{1}_{g' \geq 0} - \mathbb{1}_{g' \leq 0}) \sum_{k=0}^{2^{j_{0,n}}-1} \hat{e}_{j_{0,n}k}^2,$$

where

$$\hat{e}_{j_{0,n}k} = \sqrt{\frac{\pi}{2}} \frac{2^{j_{0,n}/2}}{\sqrt{n}} \sum_{i=1}^n \mathbb{1}_{j_{0,n}k}(i/n) |g(X_{(i-1)/n}^{(\alpha_n)}) g'(X_{(i-1)/n}^{(\alpha_n)})|^{1/2} |X_{i/n}^{(\alpha_n)} - X_{(i-1)/n}^{(\alpha_n)}|$$

and $\mathbb{1}_{g' \geq 0}$ indicates whether $x \to g'(x)$ is non-negative or not.

## 2.5. Convergence in law

We state in this section some limit theorems. In this context, it is convenient to use the notion of stable convergence in law; see Rényi [24], Aldous and Eagleson [2], Jacod and Shiryaev [17] and Jacod [14].

*Definition 1 (Stable convergence in law).* A sequence of variable $(X_n)_{n \in \mathbb{N}}$ converges stably in law to a variable $X$ ($X_n \to_{\mathcal{L}s} X$) if $X$ is defined on an appropriate extension $(\bar{\Omega}, \bar{\mathcal{F}}, \bar{\mathbb{P}})$ of $(\Omega, \mathcal{F}, \mathbb{P})$ and if for any $\mathcal{F}$-measurable bounded variable $Y$ and any bounded continuous function $g$, $\mathbb{E}[Yg(X_n)] \to \bar{\mathbb{E}}[Yg(X)]$.

For $\beta > 0$, we define the function $\Delta_\beta$ by

$$\Delta_\beta(x) = \lim_n \mathbb{E}\left[ n^{-1/2} \left( \sum_{i=1}^n Z_i \right)^2 \right],$$

with

$$Z_i = \beta(\pi/2)^{1/2} |\lfloor \{U + \beta^{-1}\sigma(x)W_{i-1}\} + \beta^{-1}\sigma(x)(W_i - W_{i-1}) \rfloor| - \sigma(x),$$

where $W$ is a Brownian motion and $U$ a uniform random variable on $[0,1]$, independent of $W$. From Delattre [8], we get that the function $\Delta_\beta$ is well defined. We have the following theorem:

**Theorem 2 (Convergence in law).** *In model (1)–(2), under Assumptions* A1*,* B *and* C1*, for $S \in \mathcal{S}$, we have the following stable convergences in law, where $B$ is a standard Brownian motion, independent of $\mathcal{F}$:*

$$\text{if } \beta_n \to 0, \qquad \sqrt{n}(\hat{\theta}_n(S) - \theta) \to_{\mathcal{L}s} \sqrt{2}(\pi - 2)^{1/2} \int_0^1 g(X_t)^2 \sigma(X_t)^2 \, \mathrm{d}B_t,$$

$$\text{if } \beta_n \to \beta > 0, \qquad \sqrt{n}(\hat{\theta}_n(S) - \theta) \to_{\mathcal{L}s} 2 \int_0^1 g(X_t)^2 \sigma(X_t) [\Delta_\beta(X_t)]^{1/2} \, \mathrm{d}B_t,$$

$$\text{if } \beta_n \to +\infty, \qquad \alpha_n^{-1}(\hat{\theta}_n(S) - \theta) \to_{\mathcal{L}s} \frac{2}{\sqrt{3}} \int_0^1 g(X_t)^2 \sigma(X_t) \, \mathrm{d}B_t.$$

## 3. Discussion

### 3.1. Comments on the results

• Our microstructure model with round-off error is obviously built to face the problem of price discreteness. Indeed, market price increments have to be multiples of the tick size. Moreover, it is striking to see how high frequency financial data do look like diffusions with round-off error; see, for example, [28]. In particular, the well-known bid-ask bounce is reproduced in this model. Furthermore, if the sampling period becomes large, the round-off error becomes insignificant. According to the theory and the empirical studies, this is also the case on the markets where it is often admitted that low frequency financial data can be modeled as data coming from a diffusion process. Hence, this model is relevant because it is clearly linked with market observations and financial theory.

• Our point of view is different from those of an additive microstructure noise. We do not make assumptions on the difference between the observed log-price and the theoretical log-price but on the observed price itself. Hence, our method is not a denoising method. We directly use the properties of the noisy data. Moreover, Li and Mykland [21] have proved that estimators built for additive noise, like the two scales estimator of Zhang, Mykland and Aït-Sahalia [31], are not robust in the case of a "quite big" rounding error.

• The estimation rates are the same as those obtained by Delattre [8] for other procedures on this model. In particular, if the order of magnitude of the round-off error is smaller than $n^{-1/2}$, we find the classical parametric rate.

• More general forms of stochastic volatility seem difficult to treat with our wavelet technique. Indeed, our proof of the central limit theorem relies on the fact that under an equivalent measure, the process can be written as a function of a Brownian motion, which is not the case for general stochastic volatility models. Nevertheless, Theorem 1 remains true in the case where the instantaneous volatility is of the form $\sigma(x,t) = g_1(x)g_2(t)$, with $g_1$ and $g_2$ as two positive functions such that $g_1 \in \mathcal{C}^2(\mathbb{R})$ and $g_2 \in \mathcal{C}^1([0,1])$; see [28].
• The integrated volatility cannot be recovered in a pure rounding framework where the sequence $\alpha_n$ is supposed to be constant as $n$ goes to infinity. Nevertheless, it can be done for some particular rounding procedures; see Jacod *et al.* [15].

## 3.2. Intuition for the results and important ideas

To give some intuition for the results, introduce important ideas used in the proofs, and explain why methods based on the quadratic variation do not work here, we recall and explain an inspiring result of Delattre [8] when $\beta_n$ tends to infinity.

### 3.2.1. The behavior of the p-variations

Let $h$ be the density of a standard Gaussian variable and

$$\gamma_p(\sigma, \beta) = \int_0^1 \mathrm{d}u \int_{\mathbb{R}} \mathrm{d}y\, h(y) |(\beta u + \sigma y)^{(\beta)}|^p.$$

It is shown in Delattre [8] that if $\beta_n$ tends to infinity, that is, if the round-off error is quite big, for $p > 0$, we have

$$\alpha_n^{-p} \beta_n n^{-1} \sum_{i=1}^n |X_{i/n}^{(\alpha_n)} - X_{(i-1)/n}^{(\alpha_n)}|^p - \beta_n^{1-p} \int_0^1 \gamma_p(\sigma(X_s), \beta_n)\,\mathrm{d}s$$

tends to zero in probability. The stable convergence in law of this sequence normalized by $\alpha_n^{-1}$ has also been proved in [8].

### 3.2.2. Remarks and explanations

The point is to remark that if $p = 1$,

$$\beta_n^{1-p} \int_0^1 \gamma_p(\sigma(X_s), \beta_n)\,\mathrm{d}s = (2/\pi)^{1/2} \int_0^1 \sigma(X_s)\,\mathrm{d}s$$

and that if $p > 0$,

$$\beta_n^{1-p} \int_0^1 \gamma_p(\sigma(X_s), \beta_n)\,\mathrm{d}s - (2/\pi)^{1/2} \int_0^1 \sigma(X_s)\,\mathrm{d}s$$

tends to zero in probability. So, in the case $\beta_n$ tends to infinity, if $p$ is not equal to one, contrary to what happens when considering non-noisy data, $\sigma(X_t)^p$ does not appear in the limit of the sum of the rescaled rounded increments to the power $p$. Thus, estimating

$$\int_0^1 \sigma^2(X_s)\,\mathrm{d}s,$$

seems more complicated than estimating

$$\int_0^1 \sigma(X_s)\,\mathrm{d}s.$$

We now give an intuition for this surprising behavior of the $p$-variations through a non-rigorous argument. Introducing several important ideas, we explain why, when $\beta_n$ tends to infinity,

$$\mathbb{E}_{\sigma(X_{(i-1)/n})}[|X_{i/n}^{(\alpha_n)} - X_{(i-1)/n}^{(\alpha_n)}|^p] \approx \alpha_n^p \beta_n^{-1}(2/\pi)^{1/2}\sigma(X_{(i-1)/n}),$$

where $\mathbb{E}_{\sigma(X_{(i-1)/n})}$ denotes the expectation conditional on $\sigma(X_{(i-1)/n})$. We define the fractional part of $X_t$ by $\{X_t\} = X_t - \lfloor X_t \rfloor$. First we have to remark that

$$X_{i/n}^{(\alpha_n)} - X_{(i-1)/n}^{(\alpha_n)} = \alpha_n \lfloor \{X_{(i-1)/n}/\alpha_n\} + (X_{i/n} - X_{(i-1)/n})/\alpha_n \rfloor. \tag{3}$$

Kosulajeff [18] and Tukey [29] have established that when $\alpha$ is small, $\{X/\alpha\}$ is almost independent of $X$ and approximately follows a uniform law on $[0,1]$. More precisely, the following result has been shown by Delattre and Jacod [9].

**Lemma 1 (The fractional part of a variable).** *Let $k$ be a function on $\mathbb{R}$, $\mathcal{C}^r$ ($r \geq 1$), integrable with integrable derivatives. Let $f$ be a function on $\mathbb{R} \times [0,1]$, $\mathcal{C}^r$ in the first variable and such that for $0 \leq l \leq r$, $M_l = \sup_x \int_0^1 |\frac{\partial^l}{\partial x^l} f(x,u)|\,\mathrm{d}u < +\infty$. Then*

$$\left|\int_{\mathbb{R}} k(x)\left[f(x,\{x/\alpha\}) - \int_0^1 f(x,u)\,\mathrm{d}u\right]\mathrm{d}x\right| \leq (2\alpha)^r \sup_{0\leq l\leq r} M_l \sup_{0\leq l\leq r} \int_{\mathbb{R}} \left|\frac{\partial^l}{\partial x^l} k(x)\right|\mathrm{d}x.$$

Thus, since

$$X_{i/n} - X_{(i-1)/n} \approx \sigma(X_{(i-1)/n})(W_{i/n} - W_{(i-1)/n}),$$

we have

$$\mathbb{E}_{\sigma(X_{(i-1)/n})}[|X_{i/n}^{(\alpha_n)} - X_{(i-1)/n}^{(\alpha_n)}|^p] \approx \alpha_n^p \mathbb{E}_{\sigma(X_{(i-1)/n})}[|\lfloor U + \beta_n^{-1}\sigma(X_{(i-1)/n})Y\rfloor|^p],$$

where $U$ is a uniform variable on $[0,1]$, independent of $X$, and $Y$ is a standard Gaussian variable, independent of $X$ and $U$. Therefore, if $\beta_n$ tends to infinity,

$$\mathbb{E}_{\sigma(X_{(i-1)/n})}[|X_{i/n}^{(\alpha_n)} - X_{(i-1)/n}^{(\alpha_n)}|^p] \approx \alpha_n^{p-1}\mathbb{E}_{\sigma(X_{(i-1)/n})}[|X_{i/n}^{(\alpha_n)} - X_{(i-1)/n}^{(\alpha_n)}|].$$

We conclude our argument using the simple but nice fact that if $U$ is a uniform variable on $[0,1]$ and $Z$ is independent of $U$, with a density with respect to the Lebesgue measure,

$$\mathbb{E}[|\lfloor U + Z \rfloor|] = \mathbb{E}[|Z|].$$

## 4. Proofs

In all the proofs we use the previously defined notation. For technical reasons, we suppose without loss of generality that for given $j$, $n2^{-j}$ is a positive integer. In the following, $c$ and $c_p$ denote constants not depending on $n$, $j$ or $k$ and that may vary from line to line.

### 4.1. Preliminaries for the proofs of Theorems 1 and 2

*4.1.1. Localization procedure*

We slightly adapt here a classical localization procedure used, for example, in Delattre [8]. It will enable us to replace Assumptions B, C and C1 with much stronger assumptions in the proofs of Theorems 1 and 2. We fix two sequences $(\nu_q)_{q\geq 1}$ and $(\mu_q)_{q\geq 1}$ such that $(\nu_q)$ is strictly decreasing to $\nu$, $(\mu_q)$ is strictly increasing to $\mu$ and $(\mu_q - \nu_q) > 0$. We also fix a sequence of functions $\chi_q : \mathbb{R} \to [0,1]$ such that $\chi_q \in \mathcal{C}_b^2(\mathbb{R})$, $(\chi_q)^{1/2} \in \mathcal{C}_b^2(\mathbb{R})$ and

$$\chi_q(x) = 1 \quad \text{on} \quad [\nu_q, \mu_q] \quad \text{and} \quad \chi_q(x) = 0 \quad \text{on} \ (-\infty, \nu_{q+1}] \cup [\mu_{q+1}, +\infty).$$

For $q \in \mathbb{N}$, we define on $\mathbb{R}$ the real functions $\sigma_q$ and $g_q$ by

$$\sigma_q(x) = \sigma(x)\chi_q(x) + (1 - \chi_q(x)), \qquad g_q(x) = g(\nu_q) + \int_{\nu_q}^x g'(s)\chi_q(s)\,\mathrm{d}s.$$

We finally set

$$T_q = \inf\left\{ t \in [0,1], X_t \leq \nu_q + \alpha_q \text{ or } X_t \geq \mu_q \text{ or } \int_0^t a_s^2 \,\mathrm{d}s \geq q \right\} \wedge 1.$$

Under Assumption B, $T_q$ tends almost surely to 1 and $\mathbb{P}(T_q = 1) \to 1$ as $q \to +\infty$. Let $(W_t^q, t \geq 0)$ be defined by $W_t^q = W_{(T_q+t)\wedge 1} - W_{T_q}$ and $(Y_t^q)_{t\geq 0}$ be the solution of

$$\mathrm{d}Y_t^q = \sigma_q(Y_t^q)\,\mathrm{d}W_t^q, \qquad Y_0^q = X_{T_q}.$$

Consider now the process $(X_t^q)_{t\in[0,1]}$ defined by $X_t^q = X_t$ for $t \in [0, T_q]$ and $X_t^q = Y_{t-T_q}^q$ for $t \in (T_q, 1]$. This process satisfies

$$\mathrm{d}X_t^q = \sigma_q(X_t^q)\,\mathrm{d}W_t + a_t^q\,\mathrm{d}t,$$

where $a_t^q = a_t$ for $t \in [0, T_q]$ and $a_t^q = 0$ for $t \in (T_q, 1]$. The process $X^q$ coincides with the initial process $X$ on $[0, T_q]$. Let $q_0 = \inf\{q, \mu_q > \nu_q + \alpha_q\}$. For $q \geq q_0$, on $[0, T_q]$, $g_q(X_t^q)$ coincides with $g(X_t)$ and for $n \geq q$, $g_q(X_t^{q(\alpha_n)})$ coincides with $g(X_t^{(\alpha_n)})$. Finally, under Assumption C1, $(g_q)'$ is of constant sign on $\mathbb{R}$, $|(g_q)'|^{1/2} \in \mathcal{C}_b^2(\mathbb{R})$ and on $[0, T_q]$, $(g_q)'(X_t^q)$ coincides with $g'(X_t)$. Furthermore, for $n \geq q$, $(g_q)'(X_t^{q(\alpha_n)})$ coincides with $g'(X_t^{(\alpha_n)})$.

Hence it is enough to prove Theorems 1 and 2 for the processes $X^q$, for all $q \geq q_0$, and so it is enough to prove Theorem 1 under Assumptions B' and C' instead of Assumptions B and C and Theorem 2 under Assumptions B' and C1' instead of Assumptions B and C1, with Assumptions B', C' and C1' defined the following way:

**Assumption B′.**

(i) There exists $c > 0$ such that for all $x \in \mathbb{R}$, $\sigma(x) \geq c$,
(ii) $x \to \sigma(x) \in \mathcal{C}_b^2(\mathbb{R})$,
(iii) $\sup_{\omega \in \Omega} \int_0^1 a_s^2 \, ds < +\infty$.

**Assumption C′.**

$x \to g(x) \in \mathcal{C}_b^2(\mathbb{R})$ and there exists $c > 0$ such that for all $x \in \mathbb{R}$, $g(x) \geq c$.

**Assumption C1′.**

(i) $x \to g(x) \in \mathcal{C}_b^2(\mathbb{R})$ and there exists $c > 0$ such that for all $x \in \mathbb{R}$, $g(x) \geq c$,
(ii) $x \to g'(x)$ is of constant sign on $\mathbb{R}$ and $x \to |g'(x)|^{1/2} \in \mathcal{C}_b^2(\mathbb{R})$.

*4.1.2. Change of probability*

Under Assumption B′, by the Girsanov theorem, we can construct a probability $\mathbb{P}'$ on $(\Omega, \mathcal{F})$, absolutely continuous with respect to $\mathbb{P}$ and a Brownian motion under $\mathbb{P}'$, $(W_t', t \geq 0)$ such that

$$dX_t = \sigma(X_t) \, dW_t' + \frac{1}{2} \sigma(X_t) \frac{\partial}{\partial x} \sigma(X_t) \, dt.$$

Assumption B′ holds for this representation. We define the following supplementary hypothesis:

**Assumption D.**

$$a_t = \frac{1}{2} \sigma(X_t) \frac{\partial}{\partial x} \sigma(X_t).$$

The convergence in probability and the stable convergence in law being preserved by an absolutely continuous change of probability, it is consequently enough to prove Theorem 1 under Assumptions A, B′, C′ and D and Theorem 2 under Assumptions A1, B′, C1′ and D. Under Assumptions B′ and D, $X_t = h(W_t)$ with $h : x \to S^{-1}(x + S(x_0))$ and

$$S : x \to \int_0^x \frac{1}{\sigma(y)} \, dy.$$

For simplicity, we suppose now that $x_0 = 0$. Note that $X$ is a homogeneous Markov process with transition densities

$$p_t(x, y) = \sigma(y)^{-1} (2\pi t)^{-1/2} \exp[-(2t)^{-1}(S(y) - S(x))^2].$$

Moreover, the following inequalities hold; see, for example, Delattre and Jacod [9].

$$\int \left| \frac{\partial^{i+j}}{\partial x^i \, \partial x^j} p_t(x, y) \right| dy \leq c t^{-(i+j)/2}, \qquad i + j \leq 2, \tag{4}$$

$$\int \left|\frac{\partial^i}{\partial x^i} q_t(x,y)\right| |y|^p \, dy \leq c_p t^{p/2}, \qquad i \leq 2, \tag{5}$$

with $q_t(x,y) = p_t(x, x+y)$. We now give the proofs of Theorems 1 and 2.

## 4.2. The behavior of the sampling functions

We give in this section a key proposition for the proofs of Theorems 1 and 2. As in Delattre [8], we consider the following assumption:

***Assumption E.*** Let $(x,u,y) \to f_n(x,u,y)$ be a sequence of real functions on $\mathbb{R} \times [0,1] \times \mathbb{R}$. The sequence $f_n$ satisfies Assumption E if the functions $f_n$ are twice continuously differentiable with respect to the first variable and if there exists $\gamma > 0$ such that for all $n \geq 1$,

(i) $|f_n(x,u,y)| \leq \gamma(1+\beta_n^2)(1+|y|^\gamma)$,
(ii) $\int_0^1 |f_n(x,u,y)| \, du \leq \gamma(1+|y|^\gamma)$,
(iii) $|\frac{\partial^i}{\partial x^i} f_n(x,u,y)| \leq \gamma(1+\beta_n^2)(1+|y|^\gamma), i=1,2$,
(iv) $\int_0^1 |\frac{\partial^i}{\partial x^i} f_n(x,u,y)| \, du \leq \gamma(1+|y|^\gamma), i=1,2$.

***Notation.*** For some sequences of real functions $x \to g_n(x)$ on $\mathbb{R}$ and $(x,u,y) \to f_n(x,u,y)$ on $\mathbb{R} \times [0,1] \times \mathbb{R}$, we define

$$V^{jk}(n, g_n) = \frac{2^{j/2}}{n} \sum_{i=1}^n \mathbb{1}_{jk}(i/n) g_n(X_{(i-1)/n})$$

and

$$V^{jk}(n, f_n) = \frac{2^{j/2}}{n} \sum_{i=1}^n \mathbb{1}_{jk}(i/n) f_n(X_{(i-1)/n}, \{X_{(i-1)/n}/\alpha_n\}, \sqrt{n}[X_{i/n} - X_{(i-1)/n}]).$$

Let $h_\sigma$ be the density of a centered Gaussian variable with variance $\sigma^2$. For a real function $(x,u,y) \to f(x,u,y)$ on $\mathbb{R} \times [0,1] \times \mathbb{R}$, we set

$$mf(x,u) = \int_\mathbb{R} h_{\sigma(x)}(y) f(x,u,y) \, dy, \qquad Mf(x) = \int_0^1 mf(x,u) \, du.$$

The following proposition is a general result on the behavior of the sampling functions.

**Proposition 1 (Behavior of the sampling functions).** *Let $(x,u,y) \to f_n(x,u,y)$ be a sequence of real functions on $\mathbb{R} \times [0,1] \times \mathbb{R}$ satisfying Assumption E. Under Assumptions A, B′ and D,*

$$\mathbb{E}\left[\left(V^{jk}(n, f_n) - 2^{j/2} \int_0^1 \mathbb{1}_{jk}(s) M f_n(X_s) \, ds\right)^2\right] \leq c r_n^2$$

*for* $0 \leq j \leq \lfloor \log_2 r_n^{-1} \rfloor$ *and* $0 \leq k \leq 2^j - 1$. *This holds for* $0 \leq j \leq \lfloor (1+\rho) \log_2 r_n^{-1} \rfloor$ *under Assumptions* A1, B' *and* D.

### 4.3. Proof of Proposition 1

In this proof, we widely use the methods and results developed by Delattre in [8]. We set $\rho$ to zero if only Assumptions A, B' and D are satisfied and write $\mathbb{E}_{\mathcal{F}_t}$ for the conditional expectation with respect to $\mathcal{F}_t$.

*4.3.1. Fundamental decomposition*

**Notation.** Let $s_{jk} = [2^{-j}nk + 1, \ldots, 2^{-j}n(k+1)]$. We use the following notation:

$$m_n f_n(x,u) = \int q_{1/n}(x,y) f_n(x,u,\sqrt{n}y)\,dy, \qquad M_n f_n(x) = \int_0^1 m_n f_n(x,u)\,du,$$

$$\bar{m}_n f_n(x) = m_n f_n(x, \{x/\alpha_n\}) - M_n f_n(x), \qquad l_i^n f_n(x) = \int p_{i/n}(x,y) \bar{m}_n f_n(y)\,dy.$$

We set

$$f_{i+1}^n = f_n(X_{i/n}, \{X_{i/n}/\alpha_n\}, \sqrt{n}[X_{(i+1)/n} - X_{i/n}]),$$

$$\eta_i^n(f_n) = f_i^n - M_n f_n(X_{(i-1)/n}),$$

$$\delta_i^n(f, l) = \sum_{z=i}^{n \wedge (i+l-1)} (\mathbb{E}_{\mathcal{F}_{i/n}}[\eta_z^n(f)] - \mathbb{E}_{\mathcal{F}_{(i-1)/n}}[\eta_z^n(f)])$$

and

$$\mathcal{M}_{jk}^n(f_n, l) = \frac{2^{j/2}}{n} \sum_{i=1}^n \mathbb{1}_{jk}(i/n) \delta_i^n(f_n, l),$$

$$H_{jk}^n(f_n, l) = \frac{2^{j/2}}{n} \sum_{i \in s_{jk}} [\bar{m}_n f_n(X_{i/n}) - \bar{m}_n f_n(X_{(i-1)/n})]$$

$$- \frac{2^{j/2}}{n} \sum_{i \in s_{jk}} l_{(n-i) \wedge (l-1)}^n f_n(X_{(i-1)/n}) \mathbb{1}_{2 \leq (n-i) \wedge (l-1)},$$

$$K_{jk}^n(f_n, l) = \frac{2^{j/2}}{n} \sum_{i \in s_{jk}} \sum_{z=1}^{(n-i-1) \wedge (l-2)} [l_z^n f_n(X_{i/n}) - l_z^n f_n(X_{(i-1)/n})].$$

Remark that for given $n$ and $z \in s_{jk}$,

$$\mathcal{M}_y^n = \frac{2^{j/2}}{n} \sum_{i=1}^y \mathbb{1}_{jk}(i/n) \delta_i^n(f_n, l)$$

is a $(\mathcal{F}_t)$-martingale in $y$. The following fundamental decomposition will be constantly used.

**Proposition 2 (Fundamental decomposition).**

$$V^{jk}(n, f_n) - 2^{j/2} \int_0^1 \mathbb{1}_{jk}(s) M f_n(X_s) \, \mathrm{d}s$$
$$= \mathcal{M}^n_{jk}(f_n, l) + V^{jk}(n, M_n f_n - M f_n)$$
$$+ V^{jk}(n, M f_n) - 2^{j/2} \int_0^1 \mathbb{1}_{jk}(s) M f_n(X_s) \, \mathrm{d}s - H^n_{jk}(f_n, l) - K^n_{jk}(f_n, l).$$

**Proof.** We have

$$\delta^n_i(f_n, l) = \eta^n_i(f_n) - M_n f_n(X_{i/n}) + M_n f_n(X_{(i-1)/n}) - \mathbb{E}_{\mathcal{F}_{(i-1)/n}}[f^n_i] + \mathbb{E}_{\mathcal{F}_{i/n}}[f^n_{i+1}]$$
$$- \mathbb{E}_{\mathcal{F}_{(i-1)/n}}[\eta^n_{i+1}(f_n)] + \sum_{z=i+2}^{n \wedge (i+l-1)} (\mathbb{E}_{\mathcal{F}_{i/n}}[\eta^n_z(f_n)] - \mathbb{E}_{\mathcal{F}_{(i-1)/n}}[\eta^n_z(f_n)]).$$

Using that

$$\mathbb{E}_{\mathcal{F}_{i/n}}[f^n_{i+1}] = \int q_{1/n}(X_{i/n}, y) f_n(X_{i/n}, \{X_{i/n}/\alpha_n\}, \sqrt{n} y) \, \mathrm{d}y,$$

we get

$$\delta^n_i(f_n, l) = \eta^n_i(f_n) + \bar{m}_n f_n(X_{i/n}) - \bar{m}_n f_n(X_{(i-1)/n}) - \mathbb{E}_{\mathcal{F}_{(i-1)/n}}[\mathbb{E}_{\mathcal{F}_{i/n}}[\eta^n_{i+1}(f_n)]]$$
$$+ \sum_{z=i+2}^{n \wedge (i+l-1)} (\mathbb{E}_{\mathcal{F}_{i/n}}[\mathbb{E}_{\mathcal{F}_{(z-1)/n}}[\eta^n_z(f_n)]] - \mathbb{E}_{\mathcal{F}_{(i-1)/n}}[\mathbb{E}_{\mathcal{F}_{(z-1)/n}}[\eta^n_z(f_n)]])$$
$$= \eta^n_i(f_n) + \bar{m}_n f_n(X_{i/n}) - \bar{m}_n f_n(X_{(i-1)/n}) - \mathbb{E}_{\mathcal{F}_{(i-1)/n}}[\bar{m}_n f_n(X_{i/n})]$$
$$+ \sum_{z=i+2}^{n \wedge (i+l-1)} (\mathbb{E}_{\mathcal{F}_{i/n}}[\bar{m}_n f_n(X_{(z-1)/n})] - \mathbb{E}_{\mathcal{F}_{(i-1)/n}}[\bar{m}_n f_n(X_{(z-1)/n})]).$$

Since

$$\mathbb{E}_{\mathcal{F}_{i/n}}[\bar{m}_n f_n(X_{(z-1)/n})] = \int p_{(z-1-i)/n}(X_{i/n}, y) \bar{m}_n f_n(y) \, \mathrm{d}y,$$

we obtain

$$\delta^n_i(f_n, l) = \eta^n_i(f_n) + \bar{m}_n f_n(X_{i/n}) - \bar{m}_n f_n(X_{(i-1)/n})$$
$$- l^n_1 f_n(X_{(i-1)/n}) + \sum_{z=2}^{(n-i) \wedge (l-1)} [l^n_{z-1} f_n(X_{i/n}) - l^n_z f_n(X_{(i-1)/n})].$$

Thus,

$$\delta_i^n(f_n) = \eta_i^n(f_n) + \bar{m}_n f_n(X_{i/n}) - \bar{m}_n f_n(X_{(i-1)/n})$$
$$- l_{(n-i)\wedge(l-1)}^n f_n(X_{(i-1)/n})\mathbb{1}_{2\leq(n-i)\wedge(l-1)}$$
$$+ \sum_{z=1}^{(n-i-1)\wedge(l-2)} [l_z^n f_n(X_{i/n}) - l_z^n f_n(X_{(i-1)/n})].$$

We finally get

$$V^{jk}(n,f_n) - V^{jk}(n, M_n f_n) = \frac{2^{j/2}}{n} \sum_{i=1}^n \mathbb{1}_{jk}(i/n)\eta_i^n(f_n)$$
$$= \mathcal{M}_{jk}^n(f_n, l) - H_{jk}^n(f_n, l) - K_{jk}^n(f_n, l). \quad \square$$

*4.3.2. Technical lemmas*

We prove here some useful lemmas. In particular, they will enable us to control the different terms of the decomposition. We begin with a usual Riemann approximation.

**Lemma 2 (Riemann approximation).** *Let $f \in \mathcal{C}_b^1$ and*

$$A_n = \frac{2^{j/2}}{n} \sum_{i=1}^n \mathbb{1}_{jk}(i/n) f(X_{i/n}) - 2^{j/2} \int_0^1 \mathbb{1}_{jk}(s) f(X_s)\,ds.$$

*Then,*

$$\mathbb{E}[A_n^2] \leq c 2^{-j} n^{-1}.$$

**Proof.** Let

$$\xi_n^i = 2^{j/2} \int_{(i-1)/n}^{i/n} [\mathbb{1}_{jk}(s) f(X_s) - \mathbb{1}_{jk}(i/n) f(X_{i/n})]\,ds.$$

We have

$$|\xi_n^i| \leq 2^{j/2} \int_{(i-1)/n}^{i/n} |f(X_s) - f(X_{i/n})|\,ds.$$

Since $f \in \mathcal{C}_b^1$, using the Burkholder–Davis–Gundy inequality, we get $\mathbb{E}[(\xi_n^i)^2] \leq c 2^j n^{-3}$. Now, $A_n = \sum_{i=1}^n \xi_n^i$ with $n 2^{-j}$ terms in the sum. Thus,

$$\mathbb{E}[A_n^2] \leq \sum_{i=1}^n \sum_{i'=1}^n (\mathbb{E}[(\xi_n^i)^2]\mathbb{E}[(\xi_n^{i'})^2])^{1/2} \leq c 2^{-j} n^{-1}. \quad \square$$

The following lemma is a consequence of Assumption E together with Lemma 1 and inequalities (4) and (5). Details can be found in Delattre [8].

**Lemma 3.**

$$|m_n f_n(x,u)| \leq c(1+\beta_n^2), \tag{6}$$

$$\int_0^1 \left|\frac{\partial^i}{\partial x^i} m_n f_n(x,u)\right| du + \left|\frac{\partial^i}{\partial x^i} M_n f_n(x)\right| \leq c, \quad 0 \leq i \leq 2, \tag{7}$$

$$|M_n f_n(x) - M f_n(x)| \leq c n^{-1/2}, \tag{8}$$

$$|l_i^n f_n(x)| \leq c\alpha_n^2 (1 + n/i). \tag{9}$$

We end this section with the following bounds for $H_{jk}^n$ and $K_{jk}^n$.

**Lemma 4.**

$$|H_{jk}^n(f_n, l)| \leq c 2^{j/2} n^{-1} + c 2^{j/2} \alpha_n^2 [1 + n 2^{-j}(l-1)^{-1} + (\log n)\mathbb{1}_{k=2^j-1}],$$

$$|K_{jk}^n(f_n, l)| \leq c 2^{j/2} \alpha_n^2 \log n.$$

**Proof.** From inequalities (6) and (9), we get

$$n|H_{jk}^n(f_n, l)| \leq c 2^{j/2}(1+\beta_n^2) + c 2^{j/2}\alpha_n^2 \sum_{i \in s_{jk}} n[(l-1)^{-1} + (n-i)^{-1}\mathbb{1}_{2 \leq (n-i)}].$$

We also have

$$n|K_{jk}^n(f_n, l)| = 2^{j/2} \sum_{z=1}^n \sum_{i \in s_{jk}} \mathbb{1}_{1 \leq z \leq (n-i-1) \wedge (l-2)} [l_z^n f_n(X_{i/n}) - l_z^n f_n(X_{(i-1)/n})]$$

$$\leq c 2^{j/2} \alpha_n^2 \sum_{z=1}^n (1 + n/z) \leq c 2^{j/2} \alpha_n^2 n \log n. \qquad \square$$

*4.3.3. End of the proof of Proposition 1*

Until the end of Section 4.4, we take $l = n$ and omit this index in the notation. We now bound the different terms of the fundamental decomposition. By Lemma 4, since $0 \leq j \leq \lfloor (1+\rho)\log_2 r_n^{-1}\rfloor$, we get

$$|H_{jk}^n(f_n) + K_{jk}^n(f_n)| \leq c 2^{j/2}(n^{-1} + \alpha_n^2 \log n) \leq c r_n.$$

Inequality (7) together with Lemma 2 on Riemann approximation give

$$\mathbb{E}\left[\left|V^{jk}(n, Mf_n) - 2^{j/2}\int_0^1 \mathbb{1}_{jk}(s) M f_n(X_s)\right|^2\right] \leq c 2^{-j} n^{-1}.$$

Inequality (8) gives

$$|V^{jk}(n, M_n f_n - M f_n)| \leq 2^{-j/2} n^{-1/2}.$$

We now turn to the approximation term $\mathcal{M}_{jk}^n(f_n)$. We have

$$\mathbb{E}[\mathcal{M}_{jk}^n(f_n)^2] = \frac{2^j}{n^2} \sum_{i \in s_{jk}} \mathbb{E}[\delta_i^n(f)^2].$$

From the results of Delattre [8], Chapters 7 and 8, we can show that

$$\mathbb{E}[\delta_i^n(f)^2] \leq c(n\alpha_n^2 + (1+\beta_n^2)(1+\alpha_n(n/i)^{1/2})).$$

Since $\sum_{i=1}^n i^{-1/2} \leq 2\sqrt{n}$, we have

$$\mathbb{E}[\mathcal{M}_{jk}^n(f_n)^2] \leq c\alpha_n^2 + c(1+\beta_n^2)(n^{-1} + 2^{j/2}\alpha_n n^{-1})$$
$$\leq c(\alpha_n^2 + n^{-1} + 2^{j/2}\alpha_n n^{-1} + 2^{j/2}\alpha_n^3). \tag{10}$$

Putting all the inequalities together, we obtain Proposition 1.

### 4.4. Proof of Theorem 1

Using the remark on the change of probability in Section 4.1.2, the following proposition implies Theorem 1.

**Proposition 3 ($L^1$ convergence, absolute integrated volatility).** *Let $\widetilde{\theta}_n$ be the estimator defined in Section 2.2. Under Assumptions A, B′, C′ and D,*

$$\mathbb{E}[|\widetilde{\theta}_n - \theta|] \leq cr_n,$$

*with c a constant not depending on n.*

For expository purposes, we first treat the case $g(x) = 1$.

*4.4.1. Proof of Proposition 3 in the case $g(x) = 1$*

We assume here that $g(x) = 1$ and set

$$f_n(x, u, y) = (\pi/2)^{1/2} \beta_n |\lfloor u + \beta_n^{-1} y \rfloor|.$$

This specification implies

$$Mf_n(X_s) = \sigma(X_s).$$

We begin with a lemma on the behavior of the wavelet coefficients. Let $c_{j_0 k}$, $d_{jk}$ and $\hat{c}_{j_0 k}$ be as defined in Section 2.2. Thanks to the vanishing moment of $\psi$, we easily get the following result:

**Lemma 5.**
$$c_{j_0 k}^2 \leq c2^{-j_0}, \qquad \mathbb{E}[d_{jk}^2] \leq c2^{-2j}.$$

Let
$$Z_j = \sum_{k=0}^{2^j - 1} \mathcal{M}_{jk}^n(f_n) c_{jk}, \qquad \tilde{Z}_j = \sum_{k=0}^{2^j - 1} K_{jk}^n(f_n) c_{jk}.$$

We have the following lemma:

**Lemma 6.** *Let $0 \leq j \leq \lfloor (1 + \rho) \log_2 r_n^{-1} \rfloor$, then*
$$\mathbb{E}[|Z_j| + |\tilde{Z}_j|] \leq cr_n.$$

**Proof.** We have $Z_j = Z_{j,1} + Z_{j,2}$ with
$$Z_{j,1} = \sum_{k=0}^{2^j - 1} 2^{j/2} \left( \int_{k/2^j}^{(k+1)/2^j} [\sigma(X_s) - \sigma(X_{k2^{-j}})] \, ds \right) \mathcal{M}_{jk}^n,$$

$$Z_{j,2} = \frac{1}{n} \sum_{k=0}^{2^j - 1} \sigma(X_{k2^{-j}}) \, ds \sum_{i \in s_{jk}} \delta_i.$$

We easily get $\mathbb{E}[|Z_{j,1}|] \leq cr_n$. For $Z_{j,2}$, we have
$$\mathbb{E}[|Z_{j,2}|^2] = \frac{1}{n^2} \sum_{k=0}^{2^j - 1} \sum_{k'=0}^{2^j - 1} \sum_{i \in s_{jk}} \sum_{i' \in s_{jk'}} \mathbb{E}[\sigma(X_{k2^{-j}}) \sigma(X_{k'2^{-j}}) \delta_i \delta_i'].$$

For $i \neq i'$, conditioning by $\mathcal{F}_{\max(i,i')-1/n}$, we get
$$\mathbb{E}[\sigma(X_{k2^{-j}}) \sigma(X_{k'2^{-j}}) \delta_i \delta_i'] = 0.$$

Therefore,
$$\mathbb{E}[|Z_{j,2}|^2] = \frac{1}{n^2} \sum_{k=0}^{2^j - 1} \mathbb{E}\left[ \sigma(X_{k2^{-j}})^2 \mathbb{E}_{\mathcal{F}_{k2^{-j}}} \left[ \sum_{i \in s_{jk}} \delta_i^2 \right] \right]$$
$$= 2^{-j} \sum_{k=0}^{2^j - 1} \mathbb{E}[\sigma(X_{k2^{-j}})^2 \mathbb{E}_{\mathcal{F}_{k2^{-j}}}[\mathcal{M}_{jk}^{n2}]] \leq cr_n^2.$$

For $\tilde{Z}_j$, recall that
$$K_{jk}^n(f_n) = \frac{2^{j/2}}{n} \sum_{i \in s_{jk}} \mathbb{1}_{jk}(i/n) \tilde{\delta}_i,$$

with

$$\widetilde{\delta}_i = \sum_{z=1}^{n-i-1} [l_z^n f_n(X_{i/n}) - l_z^n f_n(X_{(i-1)/n})]$$

and that

$$l_z^n f_n(X_{i/n}) = \mathbb{E}_{\mathcal{F}_{i/n}}[\bar{m}_n f_n(X_{(i+z)/n})].$$

The same method gives the result. $\square$

We now end the proof of Proposition 3. From Proposition 1 and equation (3), we can write

$$\hat{c}_{j_0,nk} = c_{j_0,nk} + \mathcal{M}^n_{j_0,nk}(f_n) + V^{j_0,nk}(n, M_n f_n - M f_n)$$
$$+ V^{j_0,nk}(n, M f_n) - 2^{j_0,n/2} \int_0^1 \mathbb{1}_{j_0,nk}(s) M f_n(X_s)\,\mathrm{d}s - H^n_{j_0,nk}(f_n) - K^n_{j_0,nk}(f_n)$$

and

$$\mathbb{E}[|\hat{c}_{j_0,nk} - c_{j_0,nk}|^2] \leq c r_n^2.$$

We have

$$\mathbb{E}[|\widetilde{\theta}_n - \theta|] \leq c\mathbb{E}\left[\sum_{j=j_0,n+1}^{+\infty}\sum_k d_{jk}^2 + \left|\sum_k c_{j_0,nk}\mathcal{R}_k^n\right|\right] + c(r_n + 2^{j_0,n}r_n^2),$$

where $\mathcal{R}_k^n$ is equal to

$$V^{j_0,nk}(n, M_n f_n - M f_n) + V^{j_0,nk}(n, M f_n) - 2^{j_0,n/2}\int_0^1 \mathbb{1}_{j_0,nk}(s)Mf_n(X_s)\,\mathrm{d}s - H^n_{j_0,nk}(f_n).$$

By Lemma 5, we have

$$\mathbb{E}\left[\sum_{j=j_0,n}^{+\infty}\sum_k d_{jk}^2\right] \leq c2^{-j_0,n}.$$

Moreover, using the preceding computations, it is easy to see that

$$\mathbb{E}\left[\left|\sum_k c_{j_0,nk}\mathcal{R}_k^n\right|\right] \leq c(2^{j_0,n}n^{-1} + 2^{j_0,n}\alpha_n^2 + n^{-1/2} + \alpha_n^2 \log n).$$

The result follows.

*4.4.2. Proof of Proposition 3 in the general case*

Let

$$\hat{c}^*_{j_0,n k} = \sqrt{\frac{\pi}{2}} \frac{2^{j_{0,n}/2}}{\sqrt{n}} \sum_{i=1}^{n} \mathbb{1}_{j_{0,n} k}(i/n) g(X_{(i-1)/n}) |X^{(\alpha_n)}_{i/n} - X^{(\alpha_n)}_{(i-1)/n}|.$$

We easily get the result in the same way as in the previous proof, remarking that

$$|\hat{c}_{j_{0,n} k} - \hat{c}^*_{j_{0,n} k}| \leq c \alpha_n \sqrt{\frac{\pi}{2}} \frac{2^{j_{0,n}/2}}{\sqrt{n}} \sum_{i=1}^{n} \mathbb{1}_{j_{0,n} k}(i/n) |X^{(\alpha_n)}_{i/n} - X^{(\alpha_n)}_{(i-1)/n}|.$$

Consequently, using Proposition 1, we obtain

$$\mathbb{E}[|\hat{c}_{j_{0,n} k} - c_{j_{0,n} k}|^2] \leq c \mathbb{E}[|\hat{c}_{j_{0,n} k} - \hat{c}^*_{j_{0,n} k}|^2] + c \mathbb{E}[|\hat{c}^*_{j_{0,n} k} - c_{j_{0,n} k}|^2] \leq c r_n^2$$

and

$$\sum_k c_{j_{0,n} k} (\hat{c}_{j_{0,n} k} - c_{j_{0,n} k}) = \sum_k c_{j_{0,n} k} (\hat{c}^*_{j_{0,n} k} - c_{j_{0,n} k}) + Z,$$

with $\mathbb{E}[|Z|] \leq c \alpha_n$. The result follows.

### 4.5. Proof of Theorem 2

In this proof, Assumptions A1, B' and D are in force for $\alpha_n$ and $X$. We also assume until the end of Section 4.5.2 that $g(x) = 1$.

*4.5.1. Compensator*

We have

$$\sum_k \hat{c}^2_{j_0 k} - \int_0^1 \sigma(X_s)^2 \, ds = \sum_k (\hat{c}_{j_0 k} - c_{j_0 k})^2 + 2 \sum_k c_{j_0 k} (\hat{c}_{j_0 k} - c_{j_0 k}) - \sum_{j \geq j_0} \sum_k d^2_{jk}.$$

The central limit theorems will be derived from the double product term. If, as previously, we choose $j_0$ such that $2^{j_0}$ is of order $r_n^{-1}$, re-normalized by $r_n^{-1}$, the two other terms do not tend to zero. Hence, we can either choose $2^j > r_n^{-1}$ and compensate the first term or choose $2^j < r_n^{-1}$ and compensate the last term. The first method is classical in quadratic functionals estimation. However, it seems difficult here. Indeed, a compensator of $\sum_k (\hat{c}_{jk} - c_{jk})^2$ requires an accurate enough estimation of the function $x \to \sigma(x)$. Consequently, we compensate the last term. This is unusual, but possible in our specific setting. Of course, a one-by-one estimation of the coefficients $d_{jk}$ is probably not suitable for building the compensator. This is simply because the error between the coefficient $d^2_{jk}$ and its estimation is of the same order as the error between the coefficient $c^2_{jk}$ and its estimate. That is why we use here the following scaling property of the wavelet coefficients.

**Lemma 7.** *Let*

$$Q_j = \sum_k d_{jk}^2, \qquad G(u) = \int_0^u \psi(u)\,du \quad and \quad c(\psi) = \int_0^1 G^2(u)\,du.$$

*We have*

$$\mathbb{E}\left[\left|2^j Q_j - c(\psi)\int_0^1 h'(W_u)^2\,du\right|\right] \leq c 2^{-j/2}.$$

**Proof.** We briefly give the main steps of the proof. Details can be found in Rosenbaum [27]. First we define

$$d'_{jk} = \int \psi_{jk} W_t\,dt.$$

We easily get that the $d'_{jk}$ are independent centered Gaussian variables such that

$$\mathbb{E}[d'^2_{jk}] = c(\psi) 2^{-2j}.$$

Let $\xi : [0,1] \to \mathbb{R}$ be a deterministic bounded function, vanishing outside the interval $[k 2^{-j_0}, k' 2^{-j_0}] \subset [0,1]$ and define

$$\Sigma_j(\xi) = 2^j \sum_{k=0}^{T(2^j - 1)} (2^j d'^2_{jk} - c(\psi) 2^{-j}) \xi_{k 2^{-j}}.$$

One can show that

$$\mathbb{E}[\Sigma_j(\xi)^2] \leq c \sup_t (\xi_t^2) |k' - k| 2^{j - j_0}.$$

Using a decomposition of the function $h'^2$ in a wavelet basis, these results enable us to prove that

$$\mathbb{E}\left[\left|2^j \sum_k (2^j d'^2_{jk} - c(\psi) 2^{-j}) h'(W_{k 2^{-j}})^2\right|\right] \leq c 2^{j/2}. \tag{11}$$

Since $\psi$ has a vanishing moment,

$$\int \psi_{jk}(t) h(W_t)\,dt \approx h'(W_{k 2^{-j}}) \int \psi_{jk}(t) W_t\,dt$$

and so

$$d_{jk}^2 \approx h'(W_{k 2^{-j}})^2 d'^2_{jk}.$$

We conclude using equation (11) together with a Riemann-type approximation. □

The following lemma shows that our method enables us to estimate the remaining coefficients accurately enough.

**Lemma 8.** *Let* $S = [a, (j_{1,n}), (j_{2,n})] \in \mathcal{S}$. *Then,*

$$r_n^{-1}\left[\sum_k (\hat{c}_{j_{1,n}k} - c_{j_{1,n}k})^2 - \sum_{j \geq j_{1,n}} \sum_k d_{jk}^2 + R_n(S)\right] \to 0.$$

**Proof.** We want to compensate

$$\sum_{j=j_{1,n}}^{\lfloor (1+a)\log_2 r_n^{-1} \rfloor} Q_j.$$

We know that for big enough $j$ and $j_{2,n} \leq j$, $Q_j$ is close to $2^{j_{2,n}-j}Q_{j_{2,n}}$. Therefore, we estimate the preceding quantity by

$$\sum_{j=j_{1,n}}^{\lfloor (1+a)\log_2 r_n^{-1} \rfloor} 2^{j_{2,n}-j}\hat{Q}_{j_{2,n}},$$

with

$$\hat{Q}_{j_{2,n}} = \sum_k \hat{d}_{j_{2,n}k}^2$$

for appropriate $j_{1,n}$ and $j_{2,n}$. Let

$$U_n = \sum_{j=j_{1,n}}^{\lfloor (1+a)\log_2 r_n^{-1} \rfloor} Q_j - \sum_{j=j_{1,n}}^{\lfloor (1+a)\log_2 r_n^{-1} \rfloor} 2^{j_{2,n}-j}\hat{Q}_{j_{2,n}}$$

and $Y = c(\psi)\int_0^1 h'(W_u)^2 \, du$. We have

$$U_n = \sum_{j=j_{1,n}}^{\lfloor (1+a)\log_2 r_n^{-1} \rfloor} (Q_j - 2^{-j}Y) + 2^{-j}(Y - 2^{j_{2,n}}Q_{j_{2,n}}) + 2^{j_{2,n}-j}(Q_{j_{2,n}} - \hat{Q}_{j_{2,n}}).$$

Using the same arguments as for Proposition 1, for $j_{1,n} \leq j \leq \lfloor (1+a)\log_2 r_n^{-1} \rfloor$, we get

$$\mathbb{E}[|\hat{d}_{jk} - d_{jk}|^2] \leq cr_n^2.$$

Hence, we also obtain

$$\mathbb{E}[|d_{jk}||\hat{d}_{jk} - d_{jk}|] \leq c2^{-j}r_n.$$

Consequently, we have

$$\mathbb{E}[|U_n|] \leq c(2^{-3j_{1,n}/2} + 2^{-(j_{1,n}+j_{2,n}/2)} + 2^{2j_{2,n}-j_{1,n}}r_n^2 + 2^{j_{2,n}-j_{1,n}}r_n).$$

As $a > 0$, it is clear that

$$r_n^{-1} \sum_{j > \lfloor (1+a)\log_2 r_n^{-1} \rfloor} Q_j \to 0.$$

□

*4.5.2. Limit theorems*

We prove in this section Theorem 2 in the case $g(x) = 1$. Let

$$f_n(x, u, y) = (\pi/2)^{1/2} \beta_n |\lfloor u + \beta_n^{-1} y \rfloor|,$$

$$q_i^n = \frac{1}{n} f_n(X_{(i-1)/n}, \{X_{(i-1)/n}/\alpha_n\}, \sqrt{n}[X_{i/n} - X_{(i-1)/n}])$$

and $z_i^n = q_i^n - \int_{(i-1)/n}^{i/n} M f_n(X_s)$. We begin with some intermediary lemmas.

**Lemma 9.** *Let*

$$T_1(j_n, n) = \sum_k \left[ 2^{j_n/2} \int \mathbb{1}_{j_n k}(s)[\sigma(X_s) - \sigma(X_{k2^{-j_n}})] \, ds \right] \left[ 2^{j_n/2} \sum_{i \in s_{j_n k}} z_i^n \right].$$

*If $2^{-j_n} + r_n^{-1} 2^{j_n/2}(n^{-1} + \alpha_n^2 \log n) + 2^{j_n} r_n \to 0$, then $r_n^{-1} T_1(j_n, n) \xrightarrow{\mathbb{P}} 0$.*

**Proof.** We write $\mathcal{M}_{jk}(l_n)$ for $\mathcal{M}_{jk}(f_n, l_n)$, $T_1$ for $T_1(j_n, n)$ and $j$ for $j_n$. We use the decomposition $T_1 = T_{11} + T_{12}$ with

$$T_{11} = \sum_k \left[ 2^{j/2} \int \mathbb{1}_{jk}(s)[\sigma(X_s) - \sigma(X_{k2^{-j}})] \, ds \right] \mathcal{M}_{jk}(l_n).$$

We easily get

$$\mathbb{E}[|T_{12}|] \leq c 2^{-j/2} r_n + c 2^{j/2}(n^{-1} + \alpha_n^2 \log n + \alpha_n^2 n 2^{-j}/l_n).$$

We take $l_n = \lfloor n/\log n \rfloor$ and therefore $r_n^{-1} \mathbb{E}[|T_{12}|]$ tends to zero. We set $F(x) = \sigma[h(x)]$. The term $T_{11}$ can be written as

$$T_{11} = \sum_k \sum_{i \in s_{jk}} \left[ 2^j \int \mathbb{1}_{jk}(s)(W_s - W_{k2^{-j}}) F'(W_{k2^{-j}}) \frac{\delta_i(l_n)}{n} \, ds \right] + 2^{j/2} \sum_k R_k \mathcal{M}_{jk}(l_n),$$

with $\mathbb{E}[|R_k|^2] \leq 2^{-4j}$. Following Delattre [8], Chapters 7 and 8, there exists $\tilde{\delta}_i(l_n)$ such that for $l_n = \lfloor n/\log n \rfloor$

$$\mathbb{E}[|\delta_i(l_n) - \tilde{\delta}_i(l_n)|^2] \leq cnr_n^2/\log n,$$

$$\mathbb{E}[\tilde{\delta}_i(l_n)^2] \leq c(1 + \beta_n^2)(1 + \alpha_n[n/i]^{1/2}).$$

Hence,

$$\mathbb{E}[\mathcal{M}_{jk}(l_n)^2] = \frac{2^j}{n^2} \sum_{i \in s_{jk}} \mathbb{E}[\delta_i(l_n)^2] \leq c r_n^2.$$

Consequently, we obtain that the expectation of the second term of $T_{11}$ is less than $2^{-j/2} r_n$. The first term can be written as $A_1 + A_2 + A_3 + A_4 + A_5$ with

$$A_1 = \sum_k \sum_{i \in s_{jk}} 2^j \int_{k2^{-j}}^{(i-1)/n} (W_s - W_{k2^{-j}}) F'(W_{k2^{-j}}) \frac{\delta_i(l_n)}{n} \, ds,$$

$$A_2 = \sum_k \sum_{i \in s_{jk}} 2^j \int_{(i-1)/n}^{i/n} (W_s - W_{k2^{-j}}) F'(W_{k2^{-j}}) \frac{\delta_i(l_n)}{n} \, ds,$$

$$A_3 = \sum_k \sum_{i \in s_{jk}} 2^j \int_{i/n}^{(k+1)2^{-j}} (W_s - W_{i/n}) F'(W_{k2^{-j}}) \frac{\delta_i(l_n)}{n} \, ds,$$

$$A_4 = \sum_k \sum_{i \in s_{jk}} 2^j [(k+1)2^{-j} - i/n](W_{i/n} - W_{(i-1)/n}) F'(W_{k2^{-j}}) \frac{\delta_i(l_n)}{n} \, ds,$$

$$A_5 = \sum_k \sum_{i \in s_{jk}} 2^j [(k+1)2^{-j} - i/n](W_{(i-1)/n} - W_{k2^{-j}}) F'(W_{k2^{-j}}) \frac{\delta_i(l_n)}{n} \, ds.$$

We easily get that $\mathbb{E}[A_1^2 + A_5^2] \leq c 2^{-j} r_n^2$. For $A_2$, we have

$$\mathbb{E}[|A_2|] \leq c \frac{2^{j/2}}{n} \left( \sup_i (\mathbb{E}[\delta_i(l_n)^2]) \right)^{1/2} \leq c 2^{j/2} (1/n + \alpha_n^2).$$

We now turn to $A_3$. We write here $\delta_i$ for $\delta_i(l_n)$. We easily obtain that $\mathbb{E}[A_3^2]$ is equal to

$$2^{2j} \sum_k \sum_{\substack{i \in s_{jk} \\ i' \in s_{jk}}} \int_{i/n}^{(k+1)2^{-j}} \int_{i'/n}^{(k+1)2^{-j}} \mathbb{E}\left[ F'(W_{k2^{-j}})^2 (W_s - W_{i/n})(W_{s'} - W_{i'/n}) \frac{\delta_i}{n} \frac{\delta_{i'}}{n} \right] ds \, ds'.$$

We consider the quantity

$$u_i = \mathbb{E}_{\mathcal{F}_{i/n}} \left[ F'(W_{k2^{-j}})^2 (W_s - W_{i/n})(W_{s'} - W_{i'/n}) \frac{\delta_i}{n} \frac{\delta_{i'}}{n} \right].$$

Suppose that $i \geq i'$ and $s' > i/n$. Then

$$u_i = F'(W_{k2^{-j}})^2 \frac{\delta_i}{n} \frac{\delta_{i'}}{n} \mathbb{E}_{\mathcal{F}_{i/n}}[(W_s - W_{i/n})(W_{s'} - W_{i'/n})]$$

$$= F'(W_{k2^{-j}})^2 \frac{\delta_i}{n} \frac{\delta_{i'}}{n} \mathbb{E}[(W_s - W_{i/n})(W_{s'} - W_{i/n})].$$

Suppose that $i \geq i'$ and $s' \leq i/n$. Then

$$u_i = F'(W_{k2^{-j}})^2 \frac{\delta_i}{n} \frac{\delta_{i'}}{n} \mathbb{E}_{\mathcal{F}_{i/n}}[(W_s - W_{i/n})(W_{s'} - W_{i'/n})] = 0.$$

Finally,

$$\mathbb{E}[A_3^2] \leq c2^j \sum_k \sum_{i \in s_{jk}} \int_{i/n}^{(k+1)2^{-j}} \int_{i/n}^{(k+1)2^{-j}} \mathbb{E}\left[F'(W_{k2^{-j}})^2 \left(\frac{\delta_i}{n}\right)^2\right] ds\, ds'.$$

Hence,

$$\mathbb{E}[A_3^2] \leq c2^{-j} r_n^2.$$

For $A_4$, consider the function $\zeta$ defined on $[0,1]$ by $\zeta(t) = 1 - t$ and $\zeta_{jk}(t) = 2^{j/2} \zeta(2^j x - k)$. We have

$$A_4 = \sum_k \sum_{i \in s_{jk}} 2^j[(k+1)2^{-j} - i/n]F'(W_{k2^{-j}})(W_{i/n} - W_{(i-1)/n}) \frac{\delta_i(l_n) - \tilde{\delta}_i(l_n)}{n}$$

$$+ \sum_k F'(W_{k2^{-j}}) \sum_{i \in s_{jk}} 2^{-j/2} \zeta_{jk}(i/n)(X_{i/n} - X_{(i-1)/n})\sigma(X_{(i-1)/n})^{-1} \frac{\tilde{\delta}_i(l_n)}{n} + R,$$

with $\mathbb{E}[|R|] \leq c(1/n + \alpha_n^{3/2})$. Using that the function $f_n$ verifies in our case

$$|f_n(x, u, y)| \leq c(1 + \beta_n)(1 + |y|),$$

following Delattre [8], Chapter 6, we can show that the quantity

$$\sqrt{n}(X_{i/n} - X_{(i-1)/n})\sigma(X_{(i-1)/n})^{-1} \tilde{\delta}_i(l_n)(1 + \beta_n)^{-1}$$

can be written as $g_n(X_{(i-1)/n}, \{X_{(i-1)/n}\alpha_n\}, \sqrt{n}[X_{i/n} - X_{(i-1)/n}])$. The function $g_n$ satisfies Assumption E. Therefore, since $Mg_n(x) = 0$, using the same arguments as in the proof of Proposition 1, we can prove that

$$\mathbb{E}\left[\left|\frac{1}{n}\sum_{i \in s_{jk}} \zeta_{jk}(i/n)g_n(X_{(i-1)/n}, \{X_{(i-1)/n}\alpha_n\}, \sqrt{n}[X_{i/n} - X_{(i-1)/n}])\right|^2\right] \leq cr_n^2.$$

Consequently,

$$\mathbb{E}\left[\left|\sum_k F'(W_{k2^{-j}}) \sum_{i \in s_{jk}} 2^{-j/2} \zeta_{jk}(i/n)(X_{i/n} - X_{(i-1)/n}) \frac{\tilde{\delta}_i(l_n)}{n}\, ds\right|\right] \leq c2^{j/2} r_n^2.$$

For the first term, we use that

$$\frac{1}{n^2} \sum_{i=1}^n \mathbb{E}[|\delta_i(l_n) - \tilde{\delta}_i(l_n)|^2] \leq cr_n^2/\log n$$

and finally

$$\mathbb{E}\left[\left|\sum_k \sum_{i \in s_{jk}} 2^j[(k+1)2^{-j} - i/n](W_{i/n} - W_{(i-1)/n})F'(W_{k2^{-j}})\frac{\delta_i(l_n) - \tilde{\delta}_i(l_n)}{n}\,\mathrm{d}s\right|\right]$$
$$\leq cr_n(\log n)^{-1/2}. \qquad \square$$

**Lemma 10.** *Let*

$$T_2(n) = \sum_i \int_{(i-1)/n}^{i/n} [\sigma(X_s) - \sigma(X_{(i-1)/n})]Mf_n(X_s)\,\mathrm{d}s.$$

*We have $r_n^{-1}T_2(n) \xrightarrow{\mathbb{P}} 0$.*

**Proof.** We write $T_2$ for $T_2(n)$ and set $c_\eta = (2/\pi)^{1/2}$. We have

$$T_2 = c_\eta \sum_i \int_{(i-1)/n}^{i/n} [\sigma(X_s) - \sigma(X_{(i-1)/n})]\sigma(X_s)\,\mathrm{d}s.$$

We can write

$$T_2 = c_\eta \sum_i \int_{(i-1)/n}^{i/n} [\sigma(X_s) - \sigma(X_{(i-1)/n})]^2\,\mathrm{d}s$$
$$+ c_\eta \sum_i \int_{(i-1)/n}^{i/n} [\sigma(X_s) - \sigma(X_{(i-1)/n})]\sigma(X_{(i-1)/n})\,\mathrm{d}s.$$

Itô's formula gives

$$T_2 = c_\eta \sum_i \int_{(i-1)/n}^{i/n} \mathrm{d}s\,\sigma(X_{(i-1)/n}) \int_{(i-1)/n}^{s} \sigma'(X_t)\sigma(X_t)\,\mathrm{d}W_t$$
$$+ c_\eta \sum_i \int_{(i-1)/n}^{i/n} \mathrm{d}s\,\sigma(X_{(i-1)/n}) \int_{(i-1)/n}^{s} \left(\frac{1}{2}\sigma'(X_t)^2\sigma(X_t) + \frac{1}{2}\sigma(X_t)^2\sigma''(X_t)\right)\mathrm{d}t + R,$$

with $\mathbb{E}[|R|] \leq c/n$. Finally, we obtain

$$T_2 = R + c_\eta \sum_i \int_{(i-1)/n}^{i/n} \mathrm{d}s\,\sigma(X_{(i-1)/n}) \int_{(i-1)/n}^{s} \sigma'(X_t)\sigma(X_t)\,\mathrm{d}W_t + R',$$

with $\mathbb{E}[|\tilde{R}'|] \leq c/n$. Let

$$\eta_i = \int_{(i-1)/n}^{i/n} \mathrm{d}s\,\sigma(X_{(i-1)/n}) \int_{(i-1)/n}^{s} \sigma'(X_t)\sigma(X_t)\,\mathrm{d}W_t.$$

For $i' < i$, $\mathbb{E}_{\mathcal{F}_{i'/n}}[\eta_i] = 0$. Hence, for given $n$, $M_i^n = \sum_{j=1}^i \eta_j^n$ is a martingale. Consequently,

$$\mathbb{E}[(M_i^n)^2] = \sum_{i=1}^n \mathbb{E}[\eta_i^2].$$

Since

$$\eta_i^2 \leq \frac{1}{n} \int_{(i-1)/n}^{i/n} \mathrm{d}s \left( \sigma(X_{(i-1)/n}) \int_{(i-1)/n}^s \sigma'(X_t)\sigma(X_t)\,\mathrm{d}W_t \right)^2,$$

we get

$$\mathbb{E}[(M_i^n)^2] \leq c/n^2. \qquad \square$$

**Lemma 11.** *Let*

$$T_3(j_n, n) = \sum_k \sum_{i \in s_{j_n k}} [\sigma(X_{k2^{-j_n}}) - \sigma(X_{(i-1)/n})] z_i^n.$$

*If* $2^{-j_n} + r_n^{-1} 2^{j_n/2}(n^{-1} + \alpha_n^2 \log n) + 2^{j_n} r_n \to 0$, *then* $r_n^{-1} T_3(j_n, n) \xrightarrow{\mathbb{P}} 0$.

**Proof.** We write $T_3$ for $T_3(j_n, n)$ and $j$ for $j_n$. We have

$$-T_3 = \sum_k \sum_{i \in s_{jk}} [\sigma(X_{(i-1)/n}) - \sigma(X_{k2^{-j}})] \left( \frac{\delta_i(l_n)}{n} + r_i + R_i \right),$$

with

$$r_i = \frac{1}{n}[\bar{m}_n f_n(X_{i/n}) - \bar{m}_n f_n(X_{(i-1)/n})]$$
$$- \frac{1}{n} \sum_{z=1}^{n-i-1} [l_z^n f_n(X_{i/n}) - l_z^n f_n(X_{(i-1)/n})] - \frac{1}{n} l_{n-i}^n f_n(X_{(i-1)/n}) \mathbb{1}_{i \leq n-2}$$

and $|R_i| \leq cn^{-3/2}$. We easily get

$$\mathbb{E}\left[ \left| \sum_k \sum_{i \in s_{jk}} [\sigma(X_{(i-1)/n}) - \sigma(X_{k2^{-j}})] \frac{\delta_i(l_n)}{n} \right|^2 \right]$$
$$= \mathbb{E}\left[ \left| \sum_k \sum_{i \in s_{jk}} [\sigma(X_{(i-1)/n}) - \sigma(X_{k2^{-j}})]^2 \left( \frac{\delta_i(l_n)}{n} \right)^2 \right| \right] \leq c 2^{-j} r_n^2.$$

The second term of the decomposition can be written as

$$\sum_k \sum_{i \in s_{jk}} [\sigma(X_{(i-1)/n}) - \sigma(X_{k2^{-j}})] r_i$$

$$= B_1 + B_2 + B_3$$

with

$$B_1 = -\sum_k \sum_{i \in s_{jk}} [\sigma(X_{i/n}) - \sigma(X_{(i-1)/n})] r_i,$$

$$B_2 = \sum_k [\sigma(X_{(k+1)2^{-j}}) - \sigma(X_{k2^{-j}})] \sum_{i \in s_{jk}} r_i,$$

$$B_3 = -\sum_k \sum_{i \in s_{jk}} [\sigma(X_{(k+1)2^{-j}}) - \sigma(X_i)] r_i.$$

Using Lemma 4, we obtain

$$\mathbb{E}[|B_2|] \leq 2^{j/2} \left( \frac{1}{n} + \alpha_n^2 \log n \right).$$

For $B_1$ we consider the decomposition $B_1 = B_{11} - B_{12}$ with

$$B_{11} = \sum_i [\sigma(X_{i/n}) - \sigma(X_{(i-1)/n})] z_i^n,$$

$$B_{12} = \sum_i [\sigma(X_{i/n}) - \sigma(X_{(i-1)/n})] \left( \frac{\delta_i(l_n)}{n} + R_i \right).$$

Using the same method as for $A_4$, we get $\mathbb{E}[|B_{12}|] \leq c r_n (\log n)^{-1/2}$. We have for the other term

$$B_{11} = \sum_i (X_{i/n} - X_{(i-1)/n}) \sigma'(X_{(i-1)/n}) q_i^n$$

$$+ \sum_i (X_{i/n} - X_{(i-1)/n}) \sigma'(X_{(i-1)/n}) \frac{\sigma(X_{(i-1)/n})}{n}$$

$$+ R$$

with $\mathbb{E}[|R|] \leq c/n$. Therefore, we easily get that $\mathbb{E}[|B_{11}|] \leq c r_n n^{-1/2}$. We now treat $B_3$. The quantity

$$\sum_k \sum_{i \in s_{jk}} [\sigma(X_{(k+1)2^{-j}}) - \sigma(X_{i/n})] \frac{1}{n} [\bar{m}_n f_n(X_{i/n}) - \bar{m}_n f_n(X_{(i-1)/n})]$$

can be written as

$$\sum_k \sum_{i \in s_{jk}} [\sigma(X_{(k+1)2^{-j}}) - \sigma(X_{i/n})] \frac{1}{n} \bar{m}_n f_n(X_{i/n})$$

$$-\sum_k \sum_{i \in s_{jk}} [\sigma(X_{(k+1)2^{-j}}) - \sigma(X_{i-1/n})] \frac{1}{n} \bar{m}_n f_n(X_{(i-1)/n})$$

$$+\sum_k \sum_{i \in s_{jk}} [\sigma(X_{i/n}) - \sigma(X_{(i-1)/n})] \frac{1}{n} \bar{m}_n f_n(X_{(i-1)/n}).$$

This is equal to

$$-\sum_k [\sigma(X_{(k+1)2^{-j}}) - \sigma(X_{k2^{-j}})] \frac{1}{n} \bar{m}_n f_n(X_{k2^{-j}})$$

$$+\sum_k \sum_{i \in s_{jk}} (X_{i/n} - X_{(i-1)/n}) \sigma'(X_{(i-1)/n}) \frac{1}{n} \bar{m}_n f_n(X_{(i-1)/n})$$

$$+ R$$

with $\mathbb{E}[|R|] \leq c r_n^2$. Eventually this term is equal to

$$R' + \sum_k \sum_{i \in s_{jk}} (X_{i/n} - X_{(i-1)/n}) \sigma'(X_{(i-1)/n}) \frac{1}{n} \bar{m}_n f_n(X_{(i-1)/n}) + R,$$

with $\mathbb{E}[|R'|] \leq c 2^{j/2} r_n^2$. The quantity

$$\sqrt{n}(X_{i/n} - X_{(i-1)/n}) \sigma'(X_{(i-1)/n}) \bar{m}_n f_n(X_{(i-1)/n}) (1 + \beta_n)^{-1}$$

can be written as $g_n(X_{(i-1)/n}, \{X_{(i-1)/n} \alpha_n\}, \sqrt{n}[X_{i/n} - X_{(i-1)/n}])$. This function satisfies Assumption E. Hence, since $M g_n(x) = 0$, we obtain

$$\mathbb{E}\left[ \left| \sum_k \sum_{i \in s_{jk}} (X_{i/n} - X_{(i-1)/n}) \sigma'(X_{(i-1)/n}) \frac{1}{n} \bar{m}_n f_n(X_{(i-1)/n}) \right| \right]$$

$$\leq c r_n^2.$$

We now treat

$$\sum_k \sum_{i \in s_{jk}} [\sigma(X_{(k+1)2^{-j}}) - \sigma(X_{i/n})] \frac{1}{n} \sum_{z=1}^{(n-i-1) \wedge (l_n - 2)} [l_z^n f_n(X_{i/n}) - l_z^n f_n(X_{(i-1)/n})].$$

This can be written as

$$\sum_k \sum_z \sum_{i \in s_{jk}} [\sigma(X_{(k+1)2^{-j}}) - \sigma(X_{i/n})] \frac{1}{n} l_z^n f_n(X_{i/n})$$

$$-\sum_k \sum_z \sum_{i \in s_{jk}} [\sigma(X_{(k+1)2^{-j}}) - \sigma(X_{i-1/n})] \frac{1}{n} l_z^n f_n(X_{(i-1)/n})$$

$$+\sum_k \sum_z \sum_{i \in s_{jk}} [\sigma(X_{i/n}) - \sigma(X_{(i-1)/n})] \frac{1}{n} l_z^n f_n(X_{(i-1)/n}).$$

Hence, it is equal to

$$-\sum_k \sum_z [\sigma(X_{(k+1)2^{-j}}) - \sigma(X_{k2^{-j}})] \frac{1}{n} l_z^n f_n(X_{(i-1)/n})$$

$$+\sum_k \sum_z \sum_{i \in s_{jk}} (X_{i/n} - X_{(i-1)/n}) \sigma'(X_{(i-1)/n}) \frac{1}{n} l_z^n f_n(X_{(i-1)/n}) + R,$$

with $\mathbb{E}[|R|] \leq c\alpha_n^2 \log n$. This is finally equal to

$$R' + \sum_z \sum_k \sum_{i \in s_{jk}} (X_{i/n} - X_{(i-1)/n}) \sigma'(X_{(i-1)/n}) \frac{1}{n} l_z^n f_n(X_{(i-1)/n}) + R,$$

with $\mathbb{E}[|R'|] \leq c 2^{j/2} \alpha_n^2 \log n$. The quantity

$$\sqrt{n}(X_{i/n} - X_{(i-1)/n}) \sigma'(X_{(i-1)/n}) l_z^n f_n(X_{(i-1)/n}) (\alpha_n^2 [1 + n/z])^{-1}$$

can be written as $g_n(X_{(i-1)/n}, \sqrt{n}[X_{i/n} - X_{(i-1)/n}])$. This function satisfies Assumption E. Hence, since $Mg_n(x) = 0$,

$$\mathbb{E}\left[\left|\sum_k \sum_{i \in s_{jk}} (X_{i/n} - X_{(i-1)/n}) \sigma'(X_{(i-1)/n}) \frac{1}{n} l_z^n f_n(X_{(i-1)/n})\right|\right] \leq cn^{-1} \alpha_n^2 (1 + n/z).$$

Eventually,

$$\mathbb{E}\left[\left|\sum_z \sum_k \sum_{i \in s_{jk}} (X_{i/n} - X_{(i-1)/n}) \sigma'(X_{(i-1)/n}) \frac{1}{n} l_z^n f_n(X_{(i-1)/n})\right|\right] \leq c\alpha_n^2 (1 + \log n).$$

It is also clear that

$$\mathbb{E}\left[\left|[\sigma(X_{(k+1)2^{-j}}) - \sigma(X_{i/n})] \frac{1}{n} l_{n-i}^n f_n(X_{(i-1)/n})\right|\right] \leq c\alpha_n^2 (1 + \log n). \qquad \square$$

We finally have the following result:

**Lemma 12.** *Let*

$$C_{j_n} = \sum_k c_{j_n k}(\hat{c}_{j_n k} - c_{j_n k}).$$

*If $2^{-j_n} + r_n^{-1} 2^{j_n/2}(n^{-1} + \alpha_n^2 \log n) + 2^{j_n} r_n \to 0$, then we have the following stable convergences in law, where $B$ is a standard Brownian motion, independent of $\mathcal{F}$:*

$$\text{if } \beta_n \to 0, \qquad \sqrt{n} C_{j_n} \to_{\mathcal{L}s} \frac{1}{\sqrt{2}}(\pi - 2)^{1/2} \int_0^1 \sigma(X_t)^2 \, dB_t,$$

$$\text{if } \beta_n \to \beta > 0, \qquad \sqrt{n} C_{j_n} \to_{\mathcal{L}s} \int_0^1 \sigma(X_t)[\Delta_\beta(X_t)]^{1/2} \, dB_t,$$

$$\text{if } \beta_n \to +\infty, \qquad \alpha_n^{-1} C_{j_n} \to_{\mathcal{L}s} \frac{1}{\sqrt{3}} \int_0^1 \sigma(X_t) \, dB_t.$$

**Proof.** We have

$$\sum_k c_{j_n k}(\hat{c}_{j_n k} - c_{j_n k}) = \sum_k \left[ 2^{j_n/2} \int \mathbb{1}_{j_n k}(s) \sigma(X_s) \, ds \right] \left[ 2^{j_n/2} \sum_{i \in s_{j_n k}} z_i^n \right]$$

$$= T_1(j_n, n) + T_2 + T_3(j_n, n) + T_4(n),$$

with

$$T_4(n) = \sum_i \sigma(X_{(i-1)/n}) q_i^n - \int \sigma(X_s) M f_n(X_s) \, ds.$$

We get the result by applying the results of Delattre [8], Chapter 2, to the term $T_4$ and using that $r_n^{-1}(T_1 + T_2 + T_3)$ tends to zero in probability. □

The proof of Theorem 2 follows using Lemma 8.

*4.5.3. Proof of Theorem 2 in the general case*

We give a sketch of the proof of the result in the general case. We have the following lemma:

**Lemma 13.** *If* $2^{-j_n} + r_n^{-1} 2^{j_n/2}(n^{-1} + \alpha_n^2 \log n) + 2^{j_n} r_n \to 0$, *then we have the following convergences in stable law, where $B$ is a standard Brownian motion, independent of $\mathcal{F}$:*

$$\text{if } \beta_n \to 0, \qquad \sqrt{n} C_{j_n} \to_{\mathcal{L}s} \frac{1}{\sqrt{2}} (\pi - 2)^{1/2} \int_0^1 g(X_t)^2 \sigma(X_t)^2 \, dB_t,$$

$$\text{if } \beta_n \to \beta > 0, \qquad \sqrt{n} C_{j_n} \to_{\mathcal{L}s} -\frac{\beta}{2} \int_0^1 g(X_s) g'(X_s) \sigma(X_s)^2 \, ds$$

$$+ \int_0^1 g(X_t)^2 \sigma(X_t)[\Delta_\beta(X_t)]^{1/2} \, dB_t,$$

$$\text{if } \beta_n \to +\infty, \qquad \alpha_n^{-1} C_{j_n} \to_{\mathcal{L}s} -\frac{1}{2} \int_0^1 g(X_s) g'(X_s) \sigma(X_s)^2 \, ds$$

$$+ \frac{1}{\sqrt{3}} \int_0^1 g(X_t)^2 \sigma(X_t) \, dB_t.$$

**Proof.** In this case, we have by analogy with Section 4.5.2

$$f_n(x, u, y) = (\pi/2)^{1/2} g(x - \alpha_n u) \beta_n |\lfloor u + y/\beta_n \rfloor|.$$

Let

$$\tilde{f}_n(x,u,y) = (\pi/2)^{1/2} g(x)\beta_n |\lfloor u + y/\beta_n \rfloor|.$$

We have

$$\sum_k c_{j_n k}(\hat{c}_{j_n k} - c_{j_n k}) = \sum_k c_{j_n k}(\hat{c}^*_{j_n k} - c_{j_n k}) + Z,$$

with $\mathbb{E}[|Z|] \leq c\alpha_n$. Hence, we easily get the result when $\beta_n$ tends to zero. We also have

$$\sum_k c_{j_n k}(\hat{c}_{j_n k} - c_{j_n k}) = \sum_k c_{j_n k}\left(\hat{c}_{j_n k} - 2^{j_n/2}\int \mathbb{1}_{j_n k}(s) M f_n(X_s)\,\mathrm{d}s\right)$$
$$+ \sum_k c_{j_n k}\left(2^{j_n/2}\int \mathbb{1}_{j_n k}(s)[M f_n(X_s) - M \tilde{f}_n(X_s)]\,\mathrm{d}s\right).$$

A bias is induced by the second term if $\beta_n$ does not tend to zero. Indeed,

$$\alpha_n^{-1}[M f_n(X_s) - M \tilde{f}_n(X_s)] \approx -\tfrac{1}{2} g'(X_s)\sigma(X_s).$$

Since

$$\sum_{k=0}^{2^{j_{0,n}}-1} \hat{e}^2_{j_{0,n} k} \xrightarrow{\mathbb{P}} \int_0^1 |g(X_s)g'(X_s)|\sigma(X_s)^2\,\mathrm{d}s,$$

we are able to compensate this bias. We conclude the proof of Theorem 2 using that if $U_n$ tends to $U$ in probability on $\Omega$ and $Y_n$ tends to $Y$ in stable law, then $(U_n, Y_n)$ tends to $(U, Y)$ in stable law. □

## Acknowledgements

I am grateful to Sylvain Delattre, Marc Hoffmann and Christian-Yann Robert for helpful discussions. I also thank the referee whose comments have substantially improved a former version of this paper.

Content: